\documentclass[reqno,11pt]{amsart}  

\usepackage{gensymb}
\usepackage{amsmath}
\usepackage{amssymb} 
\usepackage{mathtools}   
\usepackage{graphicx}
\usepackage{graphics}  
\usepackage{color}
\usepackage{verbatim} 
\usepackage[normalem]{ulem} 
\usepackage[mathscr]{eucal}
\usepackage{upgreek}
\usepackage{enumerate} 
\usepackage{cite}
\usepackage{amsmath}
\usepackage{amsfonts}
\usepackage{amssymb}
\usepackage{bbm}
\usepackage{cancel}
\usepackage{graphicx}
\usepackage{multicol}
\usepackage[utf8]{inputenc}
\usepackage{framed}
\usepackage{setspace}
\usepackage{amsthm}
\usepackage{hyperref}
\usepackage{caption}
\usepackage{subcaption}
\usepackage{bbm}
\usepackage{graphicx}
\usepackage{mwe}
\numberwithin{equation}{section}  
 
\usepackage[refpage,noprefix]{nomencl}
\usepackage{nomencl} 

\newtheorem{theorem}{Theorem}[section] 
\newtheorem{lemma}[theorem]{Lemma} 
\newtheorem{proposition}[theorem] {Proposition} 
 
\newtheorem{remark}[theorem]  {Remark}

\theoremstyle{definition}

 \usepackage{float}
\restylefloat{figure}

\DeclareMathAlphabet{\mathpzc}{OT1}{pzc}{m}{it}

\DeclarePairedDelimiter{\abs}{\lvert}{\rvert}

\renewcommand{\L} {\Lambda} %

\def\d{\delta} 
\newcommand{\e} {\varepsilon}

\def\l{\lambda}

\newfam\Bbbfam 
\font\tenBbb=msbm10 
\font\sevenBbb=msbm7 
\font\fiveBbb=msbm5 
\textfont\Bbbfam=\tenBbb 
\scriptfont\Bbbfam=\sevenBbb 
\scriptscriptfont\Bbbfam=\fiveBbb

\newcommand{\R}     {\mathbb{R}} 
\newcommand{\Z}     {\mathbb{Z}} 
\newcommand{\N}     {\mathbb{N}} 
\renewcommand{\P}   {\mathbb{P}}

\newcommand{\T}     {\mathbb{T}}

\def\1{{\mathchoice {1\mskip-4mu\mathrm l}      
{1\mskip-4mu\mathrm l} 
{1\mskip-4.5mu\mathrm l} {1\mskip-5mu\mathrm l}}} 
 
\def\comment#1{} 
\newtheoremstyle{thm}{2ex}{2ex}{\itshape\rmfamily}{} 
{\bfseries\rmfamily}{}{1.7ex}{} 
 
\newtheoremstyle{rem}{1.3ex}{1.3ex}{\rmfamily}{} 
{\itshape\rmfamily}{}{1.5ex}{} 
 
\newcommand{\Acal}   {{\mathcal A }}

\newcommand{\Ccal}   {{\mathcal C }} 
\newcommand{\Dcal}   {{\mathcal D }}

\newcommand{\Mcal}   {{\mathcal M }} 
 
\newcommand{\Ocal}   {{\mathcal O }}

\newcommand{\Rcal}   {{\mathcal R }} 
\newcommand{\Scal}   {{\mathcal S }}

\newcommand{\Vcal}   {{\mathcal V }} 
\newcommand{\Wcal}   {{\mathcal W }}

 \newcommand{\ex}{{\rm e}} 
 
\renewcommand{\d}{{\rm d}}

\newcommand{\Exp}{\mathscr{E}\kern-0.2mm{\operatorname{xp}}}
\newcommand{\Log}{\mathscr{L}\kern-0.2mm{\operatorname{og}}}

\renewcommand{\emptyset} {\varnothing} 
 
\newcommand{\p} {\partial}

\newcommand{\Tt}{{T_\tau}}

\newcommand\NoBlackBoxes{\global\overfullrule0pt}
\NoBlackBoxes

\renewcommand{\p}{\mathfrak{p}}

\newcommand{\Pn}{P^{\langle N \rangle}}
\newcommand{\En}{E^{\langle N \rangle}}
\newcommand{\Dt}{\Delta^\tau_N}
\newcommand{\Ene}{E^{\langle N \rangle}_{n,\e}}
\newcommand{\Pne}{P^{\langle N \rangle}_{n,\e}}
\newcommand{\Zt}{\Tt^{-1/2}\Z^2}
\newcommand{\Jfrak}{{\mathfrak{j}_o}}
\newcommand\mycom[2]{\genfrac{}{}{0pt}{}{#1}{#2}}

\newcommand{\hk}[1]{^{(#1)}}
\newenvironment{proofsect}[1] 
{\vskip0.1cm\noindent{\scshape #1.}\hskip0.5cm} 

\begin{document}

\title[\hfill On the intersections of random walks in two dimensions \hfill]
{A note on the intersections of two random walks in two dimensions}

\author[Quirin  Vogel]{Quirin  Vogel}
\address[Quirin  Vogel]{Mathematics Institute, University of Warwick, Coventry CV4 7AL, United Kingdom}
\email{Q.Vogel@warwick.ac.uk}

\thanks{}
  

 
\begin{abstract}
    In this note we prove a large deviation result for the intersection of the ranges of two independent random walks in dimension two. This complements the study of Phetpradap from 2011, where the intersection in dimension three and above was studied. 
\end{abstract}
 
 \maketitle
%
%
\section{Introduction}
The purpose of this short note is to close a gap in the literature: to provide a proof of the upper large deviation behaviour for the intersection of independent ranges of random walks in dimension two. The case $d\ge 3$ was settled in \cite{PP11} and was itself heavily based on \cite{van2004volume}, the celebrated paper in which the large deviation behaviour of the intersection volume of independent Wiener sausages was identified and proven. The rate functions for the intersection of random walks agree (up to a multiplicative constant) with the one given in \cite{van2004volume} and the proofs follow their set-up.\\
The intersection of independent ranges has been studied quite extensively in the past: it has been known for almost 70 years that $k$ random walks intersect infinitely often if and only if $k(d-2)\le dk$, see \cite{dvoretzky1951some}. In \cite{le1986proprietes}, a weak limit law for the intersections (scaling to Brownian mutual intersection local time) had been obtained. Moderate deviations at scales slower than the mean were obtained in \cite{chen2005moderate}. We refer the reader to \cite{chen2010random} for an excellent overview of those past results and their proofs.\\
It should be mentioned that the study of the intersection points can be seen as a natural continuation of the study of the range of random walks: indeed, the continuum result of the intersection volume is an expansion of the work in \cite{BBH}, where the large deviation behaviour of the volume of a single Wiener sausage was analysed. In \cite{PP11} the large deviation behaviour of the range (which is the natural lattice analogue of the volume of a Wiener sausage) of a random walk was characterised for $d\ge 3$. The case $d=2$ was settled in \cite{liu2019large}.\\ As some of the proofs in this work are quite similar to those in \cite{liu2019large}, we sometimes refer to that reference for a more detailed description. Furthermore, when an argument carries over directly from the case $d \ge 3$, we refer the reader to \cite{PP11}.\\
In the course of proving our main result, we rely heavily on the main result in \cite{liu2019large}, first-hitting time estimates given in \cite{uchiyama2011} and a KMT-type coupling in \cite{einmahl1989extensions}. The proof has three main steps: the first is an LDP for the intersection of the ranges on a torus with diverging volume. From this result, an upper and a lower bound are derived. Contrary to study of the range in \cite{liu2019large}, the removal of the torus restrictions is the most difficult part of the proof. Here, we exploit that the intersection points of the random walks exhibit the same clumping as the intersection volume of Wiener sausages in \cite{van2004volume}.
\section{Main Result and Setting}\label{sctn:MainRslt}
Take $\left(X_j\hk{1}\right)_{j=1}^\infty$ and $\left(X_j\hk{2}\right)_{j=1}^\infty$ two i.i.d. families of random variables with values in $\Z^2$. Suppose that every $X_j\hk{i}$ has mean zero and identity as covariance, for $i=1,2$. Let $H\colon [0,\infty)\to (0,\infty)$ be continuous and increasing and satisfy
\begin{equation}\label{MomentAss}
    \lim_{n\to\infty}\frac{1}{\log n}\log{H(n)}=\infty\, .
\end{equation}
We require that 
$$
E\left[H\left(\abs{X_1\hk{i}}\right)\right]<\infty\, ,
$$
for at least one such $H$.\\
Let our two random walks be defined as $S\hk{i}_n=\sum_{j=1}^nX\hk{i}_j$, for $i=1,2$ and $n\in\N$. The measure governing both random walks (starting at the origin) is denoted by $\P$. Let $J_n$ be the number of sites contained in both ranges of the two random walks, i.e. 
\begin{equation}
    J_n=\{x\in\Z^2\colon \forall i\in\{1,2\}\,\exists j\in\{1,\ldots, n\}\, \text{ with }S\hk{i}_j=x\}\, .
\end{equation}
We furthermore introduce the two relevant scales used in this text, both depending on $n$:
\begin{equation}
    \begin{split}
        \tau & = \log(n)\, , \\
        \Tt & =\frac{n}{\log(n)}=\frac{\ex^\tau}{\tau}\, .
    \end{split}
\end{equation}
Our main result is the following scaling limit.
\begin{theorem}\label{thm main}
Under the above-stated conditions on the random walks, we have that for $c>0$
\begin{equation}
    \lim_{n\to\infty}\frac{1}{\tau}\log\P\left(J_n\ge c\Tt\right)=-I_2(c)\, ,
\end{equation}
where 
\begin{equation}
    I_2(c)=\inf_{\phi\in\Theta(c)}\left[\int_{\R^2}\abs{\nabla{\phi}}^2(x)\d x\right]\, ,
\end{equation}
with 
\begin{equation}
    \Theta(c)=\Big\{ \phi\in H^1(\R^2)\colon \int_{\R^2}\phi^2(x)\d x=1,\, \int_{\R^2}\left(1-\ex^{-2\pi \phi^2(x)}\right)^2\d x  \ge c\Big\}\, .
\end{equation}
\end{theorem}
The rate function is the same (up to a multiplicative factor, compare \cite[Equation 1.9]{van2004volume}) as the one for the Wiener sausage case. Thus, we refer the reader to \cite[Theorem 3,4,6]{van2004volume} for its properties.
\begin{remark}
It is straight-forward to extend the results from Theorem \ref{thm main} to three or more random walks. The rate functions for these cases also agree with the ones for the Wiener sausage and are given in \cite[Section 1.6]{van2004volume}. We will not prove any results for three or more random walks in this work. 
\end{remark}
\section{Proof of Theorem \ref{thm main}}
\subsection{Proof of the torus LDP}
We begin the proof of Theorem \ref{thm main} by proving an LDP for the number of intersections on the torus:\\
Let $N>0$ be fixed and denote the continuum tours $[-N/2,N/2)^2$ by $\L_N$. The (rescaled) discrete torus $\Delta^\tau_N$ is defined as $\L_N\cap\Tt^{-1/2}\Z^2$.
Let $J_n$ be the number of intersections up to time $n$. Let $\Pn,\En$ be the measures governing the random walks projected (from $\Z^2$) onto $\Dt$. We implicitly use rounding for objects defined on the integers. Let $\Rcal\hk{i}_n$ be those sites visited by the $i-$th random walk up to time $n$.
\begin{proposition}\label{prop torus ldp}
$\Tt^{-1}J_{n}$ satisfies an LDP under $\Pn$ with rate $\tau$ and rate function $I_N$ where
\begin{equation}
    I_N(c)=\inf_{\phi\in \Theta_N(c)}\int_{\L_N}\abs{\nabla \phi}^2(x)\d x\, ,
\end{equation}
with
\begin{equation}
    \Theta_N(c)=\{\phi\in H^1(\L_N)\colon \int_{\L_N}\phi^2(x)\d x=1,\, \int_{\L_N}\left(1-\ex^{-2\pi \phi^2(x)}\right)^2\d x\ge c\}\, .
\end{equation}
\end{proposition}
\begin{proofsect}{\textbf{Proof of Proposition \ref{prop torus ldp}}}
The proof is quite similar to the proof of \cite[Proposition 1]{liu2019large}. We divide it into several steps.\\
\textbf{Skeleton approximation: }
define the skeleton walks
\begin{equation}
    \Scal_{n,\e}\hk{i}=\{S_{j\e \Tt}\hk{i}\}_{1\le j\le \tau/\e}\, ,
\end{equation}
for $i=1,2$. Let $\Pne,\Ene$ be the expectation given $\{\Scal_{n,\e}\hk{k}\}_{k=1,2}$, where $S_j\hk{i}$ is distributed with respect to $\Pn$. We then have the following result on exponential equivalence:
\begin{lemma}\label{lem exp equi}
It holds that for all $\delta>0$
\begin{equation}
    \lim_{\e\downarrow 0}\limsup_{n\to\infty}\frac{1}{\tau}\log\Pn\left(\abs{J_n-\Ene[J_n]}\ge \delta \Tt\right)=-\infty\, .
\end{equation}
\end{lemma}
The proof of the above lemma is a straight-forward extension of \cite[Proposition 4]{liu2019large} and we thus omit it. \\
\textbf{LDP for the skeleton walk:}
define $h(\mu|\nu)$ has the (cross) entropy between (real-valued) measures. Let $I_\e^{(2)}\colon \Mcal_1(\L_N\times \L_N)\to [0,\infty)$ be the function defined as follows
\begin{equation}\label{entropyFNCT}
    I_\e^{(2)}(\mu)=\begin{cases}
     h(\mu|\mu_1\otimes\p^{(N)}_\e) &\text{ if }\mu_1=\mu_2\, ,\\
    +\infty &\text{ otherwise}\, ,
    \end{cases}
\end{equation}
where $\p^{(N)}_\e=\p^{(N)}_\e(y-x)\d y$ is the measure induced by the Brownian transition kernel $ \left(\p^{(N)}_t(z)\right)_{t\ge 0}^{ z\in \L_N}$ on $\L_N$ ($x$ is with respect to $\mu_1$).\\
Fix $\eta>0$ and let $\Phi_\eta\colon \Mcal_1(\L_N\times \L_N)\times \Mcal_1(\L_N\times \L_N)\to [0,\infty)$ be defined as
\begin{equation}
\begin{split}
    \Phi_\eta(\mu_1,\mu_2)=\int_{\L_N}\d x &\left( 1-\exp\left[ -\eta  \int_{\L_N\times \L_N} 2\pi\phi_\e(y-x,z-x)\mu_1(\d y,\d z)\right]\right)\\
    &\!\!\!\!\!\!\times \left( 1-\exp\left[ -\eta  \int_{\L_N\times \L_N} 2\pi\phi_\e(y-x,z-x)\mu_2(\d y,\d z)\right]\right)\, ,
    \end{split}
\end{equation}
with
\begin{equation}\label{defPhie}
    \phi_\e(y,z)=\frac{\int_0^\e \d s\, \p_{s}^{(N)}(-y)\p_{\e-s}^{(N)}(z)}{\p_{\e}^{(N)}(z-y)}\, .
\end{equation}
We then have the following result.
\begin{lemma}\label{Coarse LDP lem}
$\left(\Tt^{-1}\Ene[J_n]\right)_{n\ge 0}$ satisfies an LDP with rate $\tau$ and rate function
\begin{equation}
    \begin{split}
        Y_\e(b)=\inf\Big\{\e^{-1}\left(I_\e^{(2)}(\mu_1)+I_\e^{(2)}(\mu_2)\right),\mu_1,\mu_2\in\Mcal_1(\L_N\times \L_N),\, \Phi_{\e^{-1}}(\mu_1,\mu_2)=b\Big\}\, .
    \end{split}
\end{equation}
\end{lemma}
The proof of Lemma \ref{Coarse LDP lem} is an extension of the proof of \cite[Proposition 2]{liu2019large} using the strategy from \cite{van2004volume}. The main idea is to again rewrite the intersection times as expectation with respect to empirical measures of the respective random walks. In that case this means to show that $\Ene[J_n]$ is close to $ \Phi_{1/\e}(L\hk{1}_{n,\e},L\hk{2}_{n,\e})$. Here, $L\hk{i}_{n,\e}$ are the pair empirical measures from \cite[Equation 3.5]{liu2019large}. The result then follows from the contraction principle together with Donsker-Varadhan theory.\\
The proof of Proposition \ref{prop torus ldp} now goes as follows: due to the exponential equivalence (at speed $\tau$) of $(J_n)_n$ and $\left(\Ene[J_n]\right)_n$ from Lemma \ref{lem exp equi} and the large deviation result from Lemma \ref{Coarse LDP lem}, it suffices to remove the discretization parameter $\e>0$. This is precisely the same situation as encountered in \cite[Section 3.5]{liu2019large} and thus we refer the reader to that article for details.
This concludes the proof of Proposition \ref{prop torus ldp}.\qed
\end{proofsect}
The proof of Theorem \ref{thm main} follows by proving a lower bound
\begin{equation}\label{eq low bnd}
    \liminf_{n\to\infty}\frac{1}{\tau}\log\P\left(J_n\ge c\Tt\right)\ge -I_2(c)\, ,
\end{equation}
and the corresponding upper bound
\begin{equation}\label{eq up bnd}
    \limsup_{n\to\infty}\frac{1}{\tau}\log\P\left(J_n\ge c\Tt\right)\le -I_2(c)\, ,
\end{equation}
\subsection{Proof of the lower bound}
The proof of Equation \eqref{eq low bnd} follows directly from Proposition \ref{prop torus ldp}. Indeed, the same argument was used in \cite{liu2019large}: condition on the event that the two random walks do not hit the boundary of the torus $\Dt$ up to time $n$. On that event, the total number of intersections on $\Z^2$ is bounded above by the number of intersections on the torus. However, as we let $N\to\infty$, the total cost of that conditioning vanishes and thus the lower bound follows.
\subsection{Proof of the upper bound}
\begin{figure*}[h]
        \centering
        \begin{subfigure}[b]{0.475\textwidth}
            \centering
            \includegraphics[width=\textwidth]{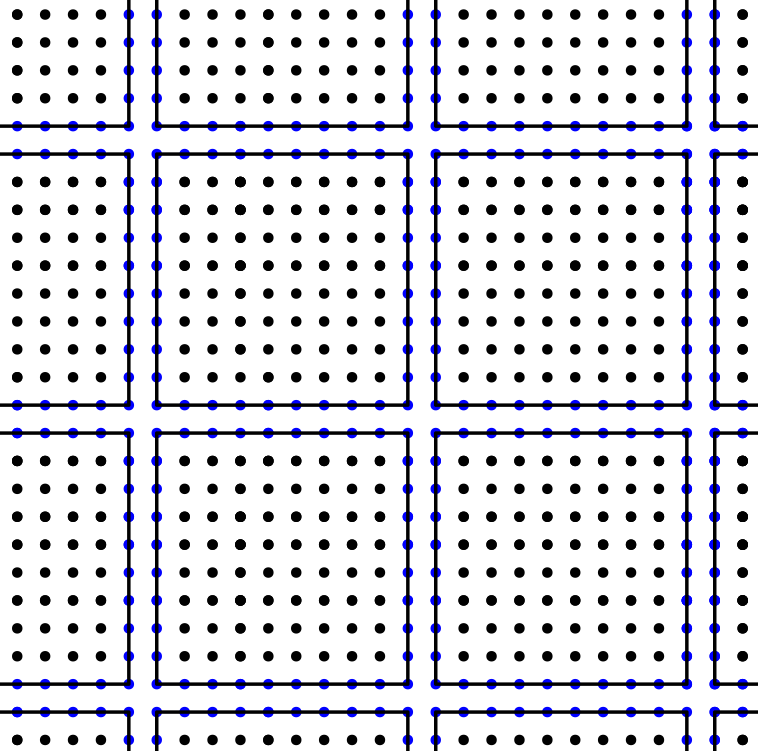}
            \label{fig:mean and std of net14}
        \end{subfigure}
        \hfill
        \begin{subfigure}[b]{0.475\textwidth}  
            \centering 
            \includegraphics[width=\textwidth]{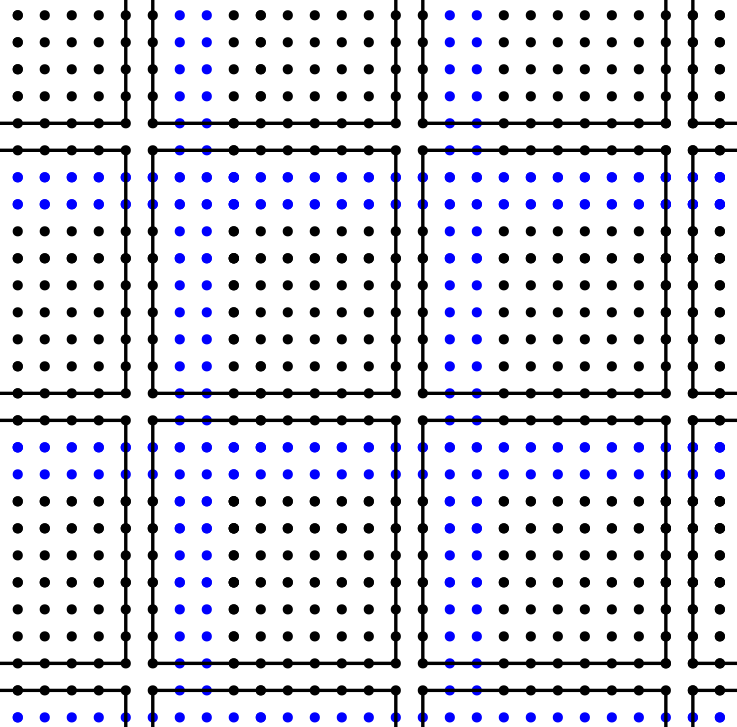}
            \label{fig:mean and std of net24}
        \end{subfigure}
        \vskip\baselineskip
        \begin{subfigure}[b]{0.475\textwidth}   
            \centering 
            \includegraphics[width=\textwidth]{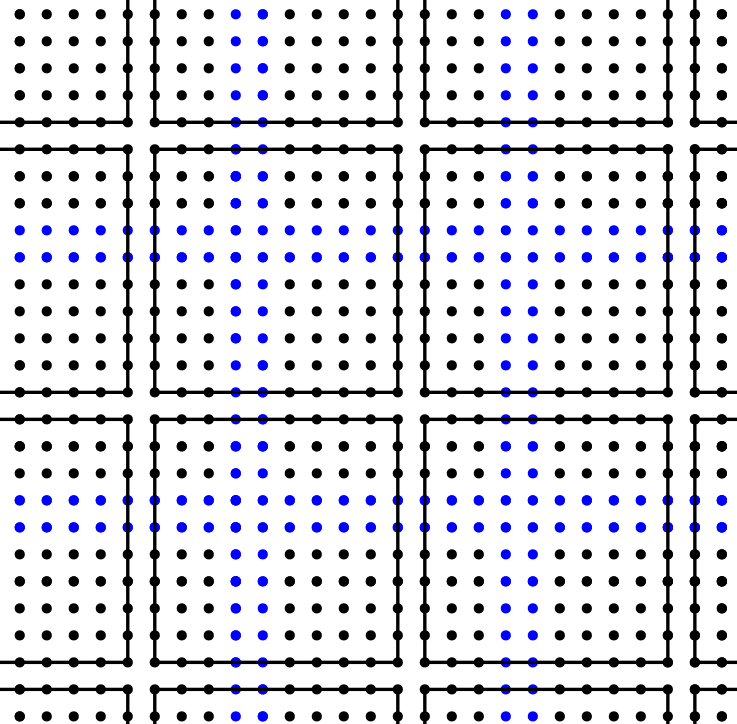}
            \label{fig:mean and std of net34}
        \end{subfigure}
        \quad
        \begin{subfigure}[b]{0.475\textwidth}   
            \centering 
            \includegraphics[width=\textwidth]{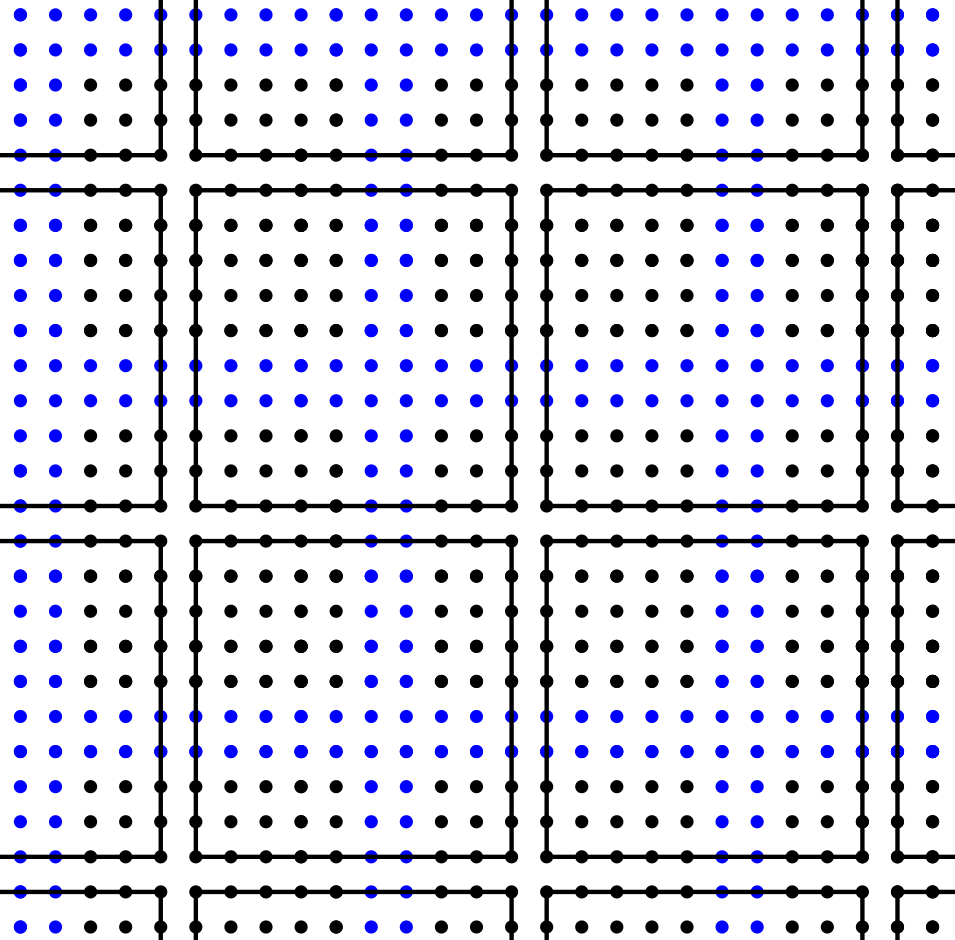}
            \label{fig:mean and std of net44}
        \end{subfigure}
        \caption[ The average and standard deviation of critical parameters ]
        {\small The different $Q_{\eta,N}^j$'s, in blue. Each blue slice has a width of $\eta\Tt^{1/2}$ points.} 
        \label{fig:IntroFigure}
    \end{figure*}

Similar to \cite{van2004volume}, we divide the proof into several steps. Firstly change notation: from now on $\Pn,\En$ are measures governing the random walks on $\Zt$. \\
\textbf{Step 1:} partition $\Tt^{-1/2}\Z^2$ into approximate $N$-boxes with each containing approximately $N^2\Tt$ many points
\begin{equation}
    \Dt(z)=\Dt+Nz\, ,
\end{equation}
for $z\in\Z^2$. For $0<\eta<N/2$, we define $Q_{\eta,N}\subset \Dt$ the $\eta/2$-neighbourhood of the faces of the boxes. We assume that $N/\eta$ is an integer. Translating $Q_{\eta,N}$ by $\eta$ in \textit{each} direction gives us $2N/\eta$ (translated) copies of $Q_{\eta,N}$ which we denote by $Q^j_{\eta,N}$, with $j=1,\ldots 2N/\eta$. Each point in $\Zt$ is contained in exactly $2$ of the $Q^j_{\eta,N}$'s. For an illustration, see Figure \ref{fig:IntroFigure}.\\
\textbf{Step 2:} define the boundary hyperplanes between the slices: for $k\in\{1,2\}$ let
\begin{equation}
    B\hk{k}=\{(z_1,z_2)\in\Zt\colon z_k=\eta\left(a+1/2\right)\text{ for }a\in \Z\}\, .
\end{equation}
Furthermore, let $B\hk{k}(x)$ be the $(z_1,z_2)\in\Zt$ such that $z_k=x$.
Define $T_1^{i,(k)}=\inf\{m>0\colon S\hk{i}_m\in  B\hk{k}\}$ and
\begin{equation}
    T_j^{i,(k)}=\inf\Big\{m>T_{j-1}^{i,(k)}\colon S\hk{i}_m\in  B_k,\, B\hk{k}\left({S\hk{i}_m}\right)\neq  B\hk{k}\left({S\hk{i}_{T_{j-1}^{i,(k)}}}\right)\Big\}\, .
\end{equation}
\textbf{Step 3:}
let 
\begin{equation}
    \Ocal_\tau=\{\forall i\in\{1,2\}\, \forall s\in\{0,\ldots,\tau\}\colon \abs{S\hk{i}_s}\le \tau^2\}\, .
\end{equation}
We then have that
\begin{equation}\label{eq ldp estm}
    \lim_{\tau \infty}\frac{1}{\tau}\log \P\left(\Ocal_\tau^c\right)=-\infty\, .
\end{equation}
Indeed, this follows similar to \cite[Proposition 1]{liu2019large}: the above statement is true when we replace the random walk with the Brownian motion. However, by the coupling from \cite[Theorem 4]{einmahl1989extensions} this can be achieved at a negligible cost at exponential scale $\tau$.\\
\textbf{Step 4:}
let $C_n\hk{k}(\eta)$ be the total number of crossings of $B\hk{k}$ made by the two random walks in direction $k$ by the time $n$ and let $C_n(\eta)=C_n\hk{1}(\eta)+C_n\hk{2}(\eta)$.
\begin{lemma}\label{lem 1}
For $M>0$
\begin{equation}
    \limsup_{\eta\to\infty}\limsup_{n\to\infty}\frac{1}{\tau}\log\Pn\left(C_n(\eta)>\frac{2M\tau}{\eta}\right)\le -C(M)\, ,
\end{equation}
with $\lim_{M\to\infty}C(M)=\infty$.
\end{lemma}
\begin{proofsect}{\textbf{Proof of Lemma \ref{lem 1}}}
Bound $\Pn\left(C_n(\eta)>\frac{2M\tau}{\eta}\right)\le 2\Pn\left(C_n\hk{1}(\eta)>\frac{2M\tau}{\eta}\right)$. As the two random walks are independent, we have
\begin{equation}
    \Pn\left(C_n\hk{1}(\eta)>\frac{2M\tau}{\eta}\right)\le 2\Pn\left(\sum_{j=1}^{M\tau (2\eta)^{-1}}T_{j+1}^{1,(k)}-T_{j}^{1,(k)}<n\right)\, .
\end{equation}
Projecting onto the first coordinate, using the reflection principle and the independence of the crossing times, we get that
\begin{equation}\label{equation combine bounds}
    \Pn\left(\sum_{j=1}^{M\tau (2\eta)^{-1}}T_{j+1}^{1,(k)}-T_{j}^{1,(k)}<n\right)\le \xi_\tau \left(\Pn\left(\max_{1\le j\le 4 \eta \Tt M^{-1}}\abs{S_j}>\eta \right)\right)^{M\tau (4\eta)^{-1}}\, ,
\end{equation}
where $\log \xi_\tau=M(2\eta)^{-1}\Tt\log 2(1+o(1))$ and we recall that $n/\tau=\Tt$. Indeed, we have that for at least half of the $j$'s that $T_{j+1}^{1,(k)}-T_{j}^{1,(k)}$ is less than $4n\eta /\left(M\tau\right)$. Denoting the combinatorial factor of choosing laf of the $j$'s by $\xi$ and applying Stirling's formula, we get the above equation. Recall that under $\Pn$ the random walk lives on the rescaled lattice $\Zt$ and thus the probability on the right-hand side converges to a finite ($\eta$ and $M$ depending) constant in $(0,1)$ by Donsker's invariance theorem. We can quantify this: due to the moment assumption in Equation \eqref{MomentAss} for the random walk, we bound
\begin{equation}\label{ldp equation bounds}
    \Pn\left(\max_{1\le j\le 4 \eta \Tt M^{-1}}\abs{S_j}>\eta \right)\le \Ocal\left(H\left(C\frac{\sqrt{M\eta}}{2}\right)\right)\, .
\end{equation}
Without loss of generality, we may choose $H$ such that
\begin{equation}
    \limsup_{n\to\infty}\frac{1}{n}\log H(n)<\infty\, .
\end{equation}
Plugging Equation \eqref{ldp equation bounds} into Equation \eqref{equation combine bounds}, taking $\log$ and letting first $n$ and then $\eta$ tend to infinity finishes the proof.\qed
\end{proofsect}
\textbf{Step 5:}
the following result follows immediately from the results in \cite{bass2009moderate} and is explicitly stated in \cite[Equation 7.2.9]{chen2010random}: for any $\l>2\pi $ one has
\begin{equation}\label{eq temp}
    \lim_{n\to\infty}\frac{1}{\tau}\log \P\left(J_n\ge \l \Tt\right)=-\infty\, .
\end{equation}
Abbreviate
\begin{equation}
    \begin{split}
        \Ccal_{\tau,M,\eta}&=\{C_\tau(\eta)\le 2M\tau \eta^{-1}\}\, ,\\
        \Vcal_\tau &=\{J_n\le 4 \pi\Tt\}\, .
    \end{split}
\end{equation}
\textbf{Step 6:}
for $j=1,\ldots,2N\eta^{-1}$ we define $C_\tau(Q^j_{\eta,N})$ the number of crossings of the hyperplanes of $Q^j_{\eta,N}$ up to time $n$. Let $J_n(Q^j_{\eta,N})$ be the number of intersection points in $Q^j_{\eta,N}$ up to time $n$. Since each points is contained in two of the $Q^j_{\eta,N}$'s, we have that on $\Ccal_{\tau,M,\eta}\cap \Vcal_\tau$ that
\begin{equation}
    \sum_{j=1}^{2N\eta^{-1}}C_\tau(Q^j_{\eta,N})\le 4 M\eta^{-1}\text{ and } \sum_{j=1}^{2N\eta^{-1}}J_n(Q^j_{\eta,N})\le 16 \pi \Tt\, .
\end{equation}
Thus, there exists a $\Jfrak \in\{1,\ldots,2N\eta^{-1}\}$ such that
\begin{equation}\label{eq spez}
    C_\tau(Q^\Jfrak_{\eta,N})\le 2d M N^{-1}\tau\,\text{ and }\, J_n(Q^\Jfrak_{\eta,N})\le 8\pi \eta N^{-1}\Tt\, .
\end{equation}
\textbf{Step 7:} choose $\eta=\sqrt{N}$ and $M=\log N$. Let $x_{N}^\Jfrak$ be the shift by which one obtains $Q^\Jfrak_{\sqrt{N},N}$ from $Q_{\sqrt{N},N}$.\\
For $z\in\Z^2$ we define
\begin{equation}
    \begin{split}
        J^\Jfrak_N(z)&=\#\{\text{ intersections in }\Delta^{\tau,\Jfrak}_N(z)\}\ ,\\
        J^\Jfrak_{N,\texttt{in}}(z)&=\#\{\text{ intersections in }Q^\Jfrak_{\sqrt{N},N}(z)\}\ ,\\
        J^\Jfrak_{N,\texttt{out}}(z)&=\#\{\text{ intersections in }\Delta^{\tau,\Jfrak}_N(z)\setminus Q^\Jfrak_{\sqrt{N},N}(z)\}\ ,
    \end{split}
\end{equation}
with $\Delta^{\tau,\Jfrak}_N(z)=\Dt+Nz+x_{N}^\Jfrak$ and $Q^\Jfrak_{\sqrt{N},N}(z)=\Dt\setminus\Delta^\tau_{N-\sqrt{N}}+Nz+x_{N}^\Jfrak$. Furthermore, let $\Rcal^{(i),\Jfrak}_\tau(z)$ be the range in $\Delta^{\tau,\Jfrak}_N(z)$ by the $i$-th walk, i.e.
\begin{equation}
    \Rcal^{(i),\Jfrak}_\tau(z)=\Rcal_n\hk{i}\cap \Delta^{\tau,\Jfrak}_N(z)\, .
\end{equation}
The highly frequented boxes $Z^\Jfrak_{\e,N}$ (identified via their shift) are defined as
\begin{equation}\label{eq highly freq boxes}
    Z^\Jfrak_{\e,N}=\{z\in\Z^2\colon \#\Rcal^{(1),\Jfrak}_\tau(z)>\e \Tt\, \text{ or }\, \#\Rcal^{(2),\Jfrak}_\tau(z)>\e \Tt\}\, .
\end{equation}
Let
\begin{equation}
    \Wcal_\tau=\{\#\Rcal_n\hk{1}\le 4\pi  \Tt,\,\#\Rcal_n\hk{2}\le 4\pi \Tt \}\, .
\end{equation}
The next two proposition imply the upper bound, this will be shown in the next step.
\begin{proposition}\label{prop 1}
There is an $N_0>0$ such that for all $0<\e<1$ and $\delta>0$ we have
\begin{equation}
    \limsup_{\tau\to \infty}\sup_{N\ge N_0}\frac{1}{\tau}\log\Pn\left(\sum_{\Z^2\setminus Z^\Jfrak_{\e,N}} J^\Jfrak_N(z)>\delta \Tt\,\text{ or }\, \sum_{Z^\Jfrak_{\e,N}} J^\Jfrak_{N,\texttt{out}}(z)>\delta\Tt\right)\le -K(\e,\delta)\, ,
\end{equation}
such that for every $\delta>0$ we have $\lim_{\e\downarrow 0}K(\e,\delta)=\infty$.
\end{proposition}
\begin{proposition}\label{prop 2}
For $N\ge 1,\, \e,\delta>0$ fixed we have
\begin{enumerate}
    \item 
        Define $N_o=2^{\#Z^\Jfrak_{\e,N}}N$. After less than $\# Z^\Jfrak_{\e,N}$ reflections in the central hyperplanes of $Q^\Jfrak_{\sqrt{N},N}$, wrapping the random walks around the torus $\Delta_{N_o}^\tau$ results in all intersections $J^\Jfrak_{N,\texttt{in}}(z)$ happen in disjoint boxes of Lebesgue measure $N^2$, $z\in Z^\Jfrak_{\e,N} $.
        \item 
        Conditioned on the event $\Ocal_\tau\cap \Ccal_{\tau,\log N,\sqrt{N}}\cap \Wcal_\tau$, the reflections have a cost (i.e. the following uniform bound on the Radon-Nikodym density) bounded by $\exp\left(\gamma_N\tau+\Ocal(\log\tau)\right)$ with $\lim_{N\to\infty}\gamma_N=0$.
\end{enumerate}
\end{proposition}
The proof of the upper bound \eqref{eq up bnd} follows directly from the two propositions above. We show that in the next step.\\
\textbf{Step 8:}
we now prove the upper bound, i.e. verify Equation \eqref{eq up bnd}. By Proposition \ref{prop 1}, Equation \eqref{eq ldp estm}, Lemma \ref{lem 1} and Equation \ref{eq temp}, we have that
\begin{equation}
    \begin{split}
        \P(J_n\ge c \Tt)\le \ex^{-\tau K(\e,\delta)}+\Pn\left(\sum_{Z^\Jfrak_{\e,N}} J^\Jfrak_{N,\texttt{out}}(z)>(c-2\delta)\Tt\cap \Ocal_\tau\cap \Ccal_{\tau,\log N,\sqrt{N}}\cap \Wcal_\tau\right)\, .
    \end{split}
\end{equation}
On account of Proposition \ref{prop 2}, we have
\begin{equation}
    \begin{split}
        \Pn &\left(\sum_{Z^\Jfrak_{\e,N}} J^\Jfrak_{N,\texttt{out}}(z)>(c-2\delta)\Tt\cap \Ocal_\tau\cap \Ccal_{\tau,\log N,\sqrt{N}}\cap \Wcal_\tau\right)\le\ex^{\gamma_N\tau+\Ocal(\log \tau)} \\
       & \times\Pn\left(\sum_{Z^\Jfrak_{\e,N}} J^\Jfrak_{N,\texttt{out}}(z)>(c-2\delta)\Tt\cap \Ocal_\tau\cap \Ccal_{\tau,\log N,\sqrt{N}}\cap \Wcal_\tau\cap \Dcal\right)\, ,
    \end{split}
\end{equation}
where $\Dcal$ is the disjointness property from Proposition \ref{prop 2}. Recall that $N_o=2^{\#Z^\Jfrak_{\e,N}}N$ and that $\#Z^\Jfrak_{\e,N}\le 4\pi\e^{-1}$ on $\Wcal_\tau$. Denote $\Pn_\T$ the measure of the random walk on the torus $\Delta_{2^{4\pi\e^{-1}}N}^\tau$. Since we are on $\Dcal$, we have that
\begin{equation}
\begin{split}
    \Pn &\left(\sum_{Z^\Jfrak_{\e,N}} J^\Jfrak_{N,\texttt{out}}(z)>(c-2\delta)\Tt\cap \Ocal_\tau\cap \Ccal_{\tau,\log N,\sqrt{N}}\cap \Wcal_\tau,\Dcal\right)\\
    &\le\Pn_\T\left(J_n>(c-2\delta)\Tt\cap \Ocal_\tau\cap \Ccal_{\tau,\log N,\sqrt{N}}\cap \Wcal_\tau\cap \Dcal\right)\, .
    \end{split}
\end{equation}
 On the torus we can apply Proposition \ref{prop torus ldp} and conclude
 \begin{equation}
     \lim_{n\to\infty}\frac{1}{\tau}\log \P\left(J_n\ge c\Tt\right)\le \max\{-K(\e,\delta),\gamma_N-I_{N_o}(c-2\delta)\}\, .
 \end{equation}
Letting $N\to\infty$, $\e\downarrow$ and then $\delta\downarrow 0$ proves Theorem \ref{thm main}. Indeed, we first use the reasoning from \cite[Section 3.6]{liu2019large} to remove the limit $N\to \infty$ in the rate function. We then use \cite[Theorem 3]{van2004volume} for the removal of the shift by $\delta$.\\
\textbf{Step 9:}
we refer the reader to \cite{PP11} for the proof of Proposition \ref{prop 2}. Indeed, the proof does not depend on the dimension and so carries over without any modifications.\\
We start the proof of Proposition \ref{prop 1}. Note that due to Equation \eqref{eq spez} and $N_0$ large enough, we can reduce the proof of Proposition \ref{prop 1} to showing that
\begin{equation}
    \limsup_{\tau\to \infty}\sup_{N\ge N_0}\frac{1}{\tau}\log\P\left(\sum_{\Z^2\setminus Z^\Jfrak_{\e,N}} J^\Jfrak_N(z)>\delta \Tt\right)\le -K(\e,\delta)\, .
\end{equation}
Define
\begin{equation}
    \Acal_{\e,N}=\{A\subset \Zt\colon \inf_{x\in\Zt}\,\,\,\sup_{z\in\Z^2}\,\#\left((A+x)\cap \Dt(z)\right)\le \e \Tt\}\, .
\end{equation}
Key for the proof of Proposition \ref{prop 1} is the next lemma.
\begin{lemma}\label{key lem}
For all $\e\in (0,1),\delta>0$
\begin{equation}
    \lim_{\e\downarrow 0}\limsup_{\tau \to \infty}\frac{1}{\tau}\log\sup_{N\ge 1}\,\sup_{A\in \Acal_{\e,N}}\Pn\left(\#(A\cap \{S_i\}_{1\le i\le n})>\delta \Tt\right)=-\infty\, ,
\end{equation}
for any $\delta>0$, where we have written $S_i$ for any of the $S_i\hk{1}$.
\end{lemma}
The proof of Proposition \ref{prop 1} follows from the above lemma in the following way: define
\begin{equation}
    A^*=\bigcup_{z\in\Z^2\colon \#(\Rcal\hk{1}_n\cap \Delta^{\tau,\Jfrak}_N(z))\le \e \Tt}\{\Rcal\hk{1}_n\cap \Delta^{\tau,\Jfrak}(z)\}\, .
\end{equation}
Note that we have by the definition of $Z^\Jfrak_{\e,N}$ in Equation \eqref{eq highly freq boxes}
\begin{equation}
\begin{split} 
    \sum_{z\in\Z^2\setminus Z^\Jfrak_{\e,N}}& J^\Jfrak_N(z)=\sum_{z\in\Z^2\setminus Z^\Jfrak_{\e,N}} \#\{\Rcal_n\hk{1}\cap\Rcal_n\hk{2}\cap\Delta_N^{\tau,\Jfrak}(z)\}\le\sum_{z\in\Z^2\setminus Z^\Jfrak_{\e,N}} \#\{\Rcal_n\hk{1}\cap\Delta_N^{\tau,\Jfrak}(z) \}\\
    &\le\sum_{z\in\Z^2} \#\{\Rcal_n\hk{1}\cap\Delta_N^{\tau,\Jfrak}(z)\cap A^* \}\le \#( A^*\cap \{S_i\}_{1\le i\le n})\, .
    \end{split}
\end{equation}
Since $A^*$ in $\Acal_{\e,N}$, Lemma \ref{key lem} implies Proposition \ref{prop 1}.\\
\textbf{Step 10:}
the proof of Lemma \ref{key lem} can be reduced further to showing
\begin{equation}
     \lim_{\e\downarrow 0}\limsup_{\tau \to \infty}\frac{1}{\tau}\log\sup_{N\ge 1}\sup_{A\in \Acal_{\e,N}}\En\left[\exp\left(\e^{-1/6}\tau\Tt^{-1}\#(A\cap \{S_i\}_{1\le i\le n})\right)\right]=0\, .
\end{equation}
Indeed, this is the exponential Chebyshev's inequality.\\
Using the subadditivity property of the range, we partition $n$ into pieces of length $\Tt$. So we bound
\begin{equation}
    \begin{split}
        \sup_{A\in \Acal_{\e,N}} & \En\left[\exp\left(\e^{-1/6}\tau\Tt^{-1}\#(A\cap \{S_i\}_{1\le i\le n})\right)\right]\\
        &\le  \sup_{A\in \Acal_{\e,N}} \En\left[\exp\left(\e^{-1/6}\tau\Tt^{-1}\sum_{k=1}^{\e^{-1/2}\tau}\#(A\cap \{S_i\}_{(k-1)\e^{1/2}\Tt\le i\le k\e ^{1/2}\Tt})\right)\right]\\
        &\le\left(  \sup_{A\in \Acal_{\e,N}}\sup_{x\in\Zt} \En_x\left[\exp\left(\e^{-1/6}\tau\Tt^{-1}\#(A\cap \{S_i\}_{1\le i\le \e ^{1/2}\Tt})\right)\right]\right)^{\e^{-1/2}\tau}\, .
    \end{split}
\end{equation}
\textbf{Step 11:}
use the inequality $\ex^u\le 1+u+\frac{u^2\ex^u}{2}$, the bound $\#(A\cap \{S_i\}_{1\le i\le \e ^{1/2}\Tt})\le \#\Rcal_{\e^{1/2}\Tt}$ and the Cauchy-Schwarz inequality to obtain
\begin{equation}\label{eq split}
\begin{split}
    \En_x &\left[\exp\left(\e^{-1/6}\tau\Tt^{-1}\#(A\cap \{S_i\}_{1\le i\le \e ^{1/2}\Tt})\right)\right]\le 1+ \e^{1/3}\tau\Tt^{-1}\En_x\left[\#(A\cap \{S_i\}_{1\le i\le \e ^{1/2}\Tt})\right]\\
    &+\frac{1}{2}\e^{2/3}\sqrt{\frac{\En_x\left[\left(\#\Rcal_{\e^{1/2}\Tt}\right)^4\right]}{\left(\e^{1/2}\Tt\tau^{-1}\right)^4}}\sqrt{\En_x\left[\exp\left(2\e^{1/6}\left(\e^{1/2}\Tt\tau^{-1}\right)^{-1}\#\Rcal_{\e^{1/2}\Tt}\right)\right]}\, .
\end{split}
\end{equation}
It follows from \cite[Theorem 6.3.1]{chen2010random} that the terms under the square roots are uniformly bounded as $n\to\infty$ (uniformly in $\e$). It remains to analyse the first term.\\
\textbf{Step 12:} expand
\begin{equation}
    \begin{split}
        \En_x&\left[\#(A\cap \{S_i\}_{1\le i\le \e ^{1/2}\Tt})\right]\le \sum_{z\in\Z^2}\En_x\left[\#\left(A\cap \{S_i\}_{1\le i\le \e ^{1/2}\Tt}\cap\Dt(z)\right)\right]\\
        &\le \sum_{z\in\Z^2}\Pn_x\left(\{S_i\}_{1\le i\le \e ^{1/2}\Tt}\cap\Dt(z)\neq \emptyset\right)\\
        &\qquad\quad\times\En_x\left[\#\left(A\cap \{S_i\}_{1\le i\le \e ^{1/2}\Tt}\cap\Dt(z)\right)\Big| \{S_i\}_{1\le i\le \e ^{1/2}\Tt}\cap\Dt(z)\neq \emptyset\right]\, .
    \end{split}
\end{equation}
Using the Markov property and spatial homogeneity, we can bound the expectation above by
\begin{equation}\label{Eq above}
    \sup_{\mycom{A\subset\Zt}{ \#(A\cap \Dt)\le \e \Tt}}\,\,\,\sup_{x\in\Zt}\,\,\,\sum_{y\in\Zt\cap A\cap\Dt}\Pn_x\left(H_y<\e^{1/2}\Tt\right)\, ,
\end{equation}
where $H_y$ is the hitting time of the point $y$. Note that by \cite[Theorem 1.7]{uchiyama2011} we can express for $\e>0$ fixed and $\abs{x}\le \e \sqrt{n}$ uniformly
\begin{equation}
    \P_0\left(H_x<n\right)=\frac{\left(1+o(1)\right) }{\log n}\int_{\abs{x}^2/(2n)}^\infty\frac{\ex^{-u}}{u}\d u=\frac{2\pi\left(1+o(1)\right) }{\log n}\int_{0}^1\p_u(\abs{x}^2/n)\d u\, .
\end{equation}
Here, $\p_t(x)$ denotes the kernel of a standard Brownian motion in $\R^d$. Since the above is decreasing in $\abs{x}$, we can bound Equation \eqref{Eq above} by
\begin{equation}
    \sum_{\mycom{x\in \Z^2}{\abs{x}\le \e \Tt^{1/2}}}P\left(H_x\le \e^{1/2}\Tt\right)\le  \frac{C_1}{\tau}\,\sum_{\mycom{x\in \Z^2}{\abs{x}\le\e\Tt^{1/2}}}\int_{0}^1\p_u\left(\abs{x}^2/(\e^{1/2}\Tt)\right)\d u\, .
\end{equation}
Approximating the sum by an integral and after a change of variables, we get
\begin{equation}
    \sum_{\mycom{x\in \Z^2}{\abs{x}\le\e\Tt^{1/2}}}\int_{0}^1\p_u\left(\abs{x}^2/(\e^{1/2}\Tt)\right)\d u\le C\e^2\Tt\int_0^{\e^{-3/2}}u^{-1}\int_0^1r \ex^{-r^2/(2u)}\d r\,\d u\, .
\end{equation}
Thus, it follows that we can bound the second term on the right-hand side of Equation \eqref{eq split} by
\begin{equation}
    \e^{1/3}\frac{\tau}{\Tt}\En_x\left[\# \left(A\cap\{S_i\}_{1\le i\le \e ^{1/2}\Tt}\right)\right]]\le C_2\e ^{5/6}\log(\e^{-1})\, .
\end{equation}
Substituting the above into Equation \eqref{eq split} concludes the proof of Lemma \ref{key lem}.\qed\\ 
We have proven the upper bound in Equation \eqref{eq up bnd} and thus have finished the proof of Theorem \ref{thm main}.\qed
\section{Acknowledgements}
The author would like to thank the Great Britain Sasakawa Foundation and the Kyushu University who funded a research visit to Japan during which some of the ideas presented above were developed. 
\bibliography{intersectionstwod}{}
\bibliographystyle{alpha}
\end{document}